\documentclass{amsart}
\usepackage{amssymb}
\usepackage{amsmath}
\usepackage{amstext}
\usepackage{amsfonts}
\usepackage{amsthm}
\usepackage{mathrsfs}
\usepackage{centernot}
\usepackage{enumitem}
\usepackage{graphicx}

\newtheorem*{theorem*}{Theorem}

\begin{document}
\title[Prime geodesic theorem]{Prime geodesic theorem of Gallagher type}
\author{Muharem Avdispahi\'{c}}
\address{University of Sarajevo, Department of Mathematics, Zmaja od Bosne
33-35, 71000 Sarajevo, Bosnia and Herzegovina}
\email{mavdispa@pmf.unsa.ba}
\subjclass[2010]{11M36, 11F72, 58J50}
\keywords{Prime geodesic theorem, Selberg zeta function, hyperbolic manifolds}

\begin{abstract}
We reduce the exponent in the error term of the prime geodesic theorem for
compact Riemann surfaces from $\frac{3}{4}$ to $\frac{7}{10}$ outside a set
of finite logarithmic measure.
\end{abstract}

\maketitle

\section{Introduction}

Let $\Gamma \subset PSL\left( 2,%
%TCIMACRO{\U{211d} }%
%BeginExpansion
\mathbb{R}
%EndExpansion
\right) $ be a strictly hyperbolic Fuchsian group and $\mathcal{F}=\Gamma
\setminus \mathcal{H}$ be the corresponding compact Riemann surface of a
genus $g\geq 2$, where $\mathcal{H}=\{z=x+iy:y>0\}$ denotes the upper
half-plane equipped with the hyperbolic metric $ds^{2}=\frac{dx^{2}+dy^{2}}{%
y^{2}}$. The Selberg zeta function on $\mathcal{F}$ is
defined by the Euler product
\begin{equation*}
Z(s)=Z_{\Gamma }(s)=\underset{\{P_{0}\}}{\prod }\underset{k=0}{\overset{\infty }%
{\prod }}(1-N(P_{0})^{-s-k})\text{, }\text{Re}(s)>1\text{,}
\end{equation*}%
where the product is taken over all primitive conjugacy classes $\left\{
P_{0}\right\} $ \ and $N(P_{0})$ is the norm of a conjugacy class $P_{0}$.
(See any of the standard references \cite{H}, \cite{I1} for a necessary
background.)

The Selberg zeta function can be continued to the whole complex plane as a
meromorphic function of a finite order. Its zeros $\rho = \frac{1}{2}%
+i\gamma $ are denumerable and closely related to the eigenvalues $\lambda$ of the
Laplace-Beltrami operator on $\mathcal{F}$. This operator
being essentially self-adjoint, its eigenvalues are non-negative and tend to
infinity. Therefore, there are finitely many of them that are less than $%
\frac{1}{4}$. We have $\gamma=\pm \sqrt{\lambda - \frac{1}{4}}$ for $\lambda
\geq \frac{1}{4}$ and $\gamma=\mp i \sqrt{\frac{1}{4}-\lambda}$ for $\lambda
< \frac{1}{4}$. So, the zeros of $Z$ are split into two parts: those lying
on the critical line $\text{Re}s=\frac{1}{2}$ and the real zeros in the
interval $[0, 1]$.

Compared to the Riemann zeta case, there are "too many zeros" of $Z$ in some
sense \cite[(6.14) on p. 113]{H}. Let $N(t)$ denote the number of zeros $\rho =\frac{1}{2}%
+i\gamma $ such that $0<\gamma \leq t$. The function $R(t)$, given by $N(t)=%
\frac{|\mathcal{F}|}{4\pi }t^{2}+R(t)$, grows as $O\left( t(\log
t)^{-1}\right) $, where $|\mathcal{F}|$ is the volume of $\mathcal{F}$.

The norm $N(P_{0})$ is determined by the length of the geodesic joining two
fixed points, necessarily the same ones for all representatives of $P_{0}$.
The statement about the number $\pi _{0}(x)$ of classes $\left\{
P_{0}\right\} $ such that $N(P_{0})\leq x$, for $x>0$, is known as the prime
geodesic theorem. It has been proved by Selberg \cite{S} and Huber \cite{H1, H2}
and subsequently generalized to various settings. The references \cite{M, D, GW, PP, D1, PS} form an interesting range of samples.

One should consult \cite[Discussion 6.20., pp. 113-115 and Discussion 15.16., pp. 253-255]{H} on the difficulties in improving Huber's $O\left( x^{\frac{3}{4}}(\log x)^{-\frac{1}{2}}\right)$. The best estimate for the error term in the prime geodesic theorem up to now
is $O\left( x^{\frac{3}{4}}(\log x)^{-1}\right)$ obtained by Randol \cite{R}. Its analogue has been also established for higher dimensional hyperbolic manifolds \cite{AG}, improving Park's theorem \cite[Th. 1.2.]{P}. An important ingredient, implicitly \cite{R} or explicitly \cite{P}, is the growth rate of the log-derivative of the Ruelle zeta vs. the Selberg zeta log-derivative \cite{AS}.

Though the expected exponent on Riemann surfaces is $\frac{1}{2}+\varepsilon$, the above
mentioned $\frac{3}{4}$ was successfully reduced only in the case of modular
surfaces $\Gamma \setminus \mathcal{H}$, $\Gamma \subset PSL(2, \mathbb{Z}) $%
. Iwaniec \cite{I} obtained $\frac{35}{48}+\varepsilon$, Luo and Sarnak \cite{LS} $%
\frac{7}{10}+\varepsilon$, Cai \cite{C} $\frac{71}{102}+\varepsilon$,
Soundararajan and Young \cite{SY} $\frac{25}{36}+\varepsilon$.

Following Gallagher's \cite{G} approach to the Riemann zeta, we shall prove
that $\frac{7}{10}$ can be achieved for $\Gamma \subset PSL(2, \mathbb{R}) $
outside a set of finite logarithmic measure.

\section{Main result}

\bigskip

Let $\psi \left( x\right) =\underset{N(P)\leq x}{\sum }\Lambda \left(
P\right) =\underset{N(P)\leq x}{\sum }\frac{\log N(P_{0})}{1-N(P)^{-1}}$ and
$\psi _{1}\left( x\right) =\int_{1}^{x}\psi \left( t\right) dt$ be the
Chebyshev resp. integrated Chebyshev function, as usual. Recall that the
error term $O\left( x^{\frac{3}{4}}(\log x)^{-1}\right) $ in the prime
geodesic theorem corresponds to $O\left( x^{\frac{3}{4}}\right) $ in the
explicit formula for $\psi $.

Our main result is given by the following theorem that substantially improves \cite[Th. 1.]{K}. and \cite[Th. 2.]{A}

\begin{theorem*}
Let $\Gamma \subset PSL(2,\mathbb{R})$ be a strictly hyperbolic
Fuchsian group. There exists a set $G $ of finite logarithmic measure such that
\begin{equation*}
\psi \left( x\right) =x+\sum_{\frac{7}{10}<\rho <1}\frac{x^{\rho }}{\rho }%
+O\left( x^{\frac{7}{10}}(\log x)^{\frac{1}{5}}\left( \log \log x\right) ^{%
\frac{1}{5}+\varepsilon }\right) \ \ \left( x\rightarrow \infty ,x\notin
G\right) \text{,}
\end{equation*}
where $\varepsilon >0$ is arbitrarily small.

\begin{proof}
As the starting point, we shall take Hejhal's explicit formula for $\psi
_{1}$ with an error term \cite[Th. 6.16. on p. 110]{H}
\begin{equation*}
\begin{split}
\psi _{1}\left( x\right) & =\alpha _{0}x+\beta _{0}x\log x+\alpha _{1}+\beta
_{1}\log x+F\left( \frac{1}{x}\right) \\
& +\frac{x^{2}}{2}+\sum_{\substack{ \rho  \\ \left\vert \gamma \right\vert
<T }}\frac{x^{\rho +1}}{\rho \left( \rho +1\right) }+O\left( \frac{x^{2}\log
x}{T}\right) \ \left( x\rightarrow \infty \right) \text{,}
\end{split}%
\end{equation*}%
where $F\left( x\right) =\left( 2g-2\right) \overset{\infty }{\underset{k=2}{%
\sum }}\frac{2k+1}{k\left( k-1\right) }x^{1-k}$.

The asymptotics of $\psi $ is conveniently derived from the asymptotics of $%
\psi _{1}$ via the relation%
\begin{equation*}
\underset{x-h}{\overset{x}{\int }}f(t)dt\leq f(x)h\leq \underset{x}{\overset{%
x+h}{\int }}f(t)dt
\end{equation*}
valid for any non-decreasing function $f$, where $h>0$.

We have%
\begin{equation}
\begin{split}
\psi \left( x\right) & \leq \frac{1}{h}\underset{x}{\overset{x+h}{\int }}%
\psi \left( t\right) dt=x+\sum_{\frac{1}{2}<\rho <1}\frac{x^{\rho }}{\rho }%
+O\left( \log x\right) +O\left( h\right) \\
& +\frac{1}{h}\left\vert \sum
_{\substack{ \text{Re}(\rho )=\frac{1}{2}  \\ \left\vert \gamma \right\vert
\leq T}}\frac{\left( x+h\right) ^{\rho +1}-x^{\rho +1}}{\rho \left( \rho
+1\right) }\right\vert +O\left( \frac{x^{2}\log x}{hT}\right)\text{.}
\end{split}
\label{f1}
\end{equation}%
Now, for $Y<T$,
\begin{equation*}
\sum_{\substack{ \text{Re}(\rho )=\frac{1}{2}  \\ \left\vert \gamma
\right\vert \leq T}}\frac{\left( x+h\right) ^{\rho +1}-x^{\rho +1}}{\rho
\left( \rho +1\right) }=\sum_{\substack{ \text{Re}(\rho )=\frac{1}{2}  \\ %
\left\vert \gamma \right\vert \leq Y}}\frac{\left( x+h\right) ^{\rho
+1}-x^{\rho +1}}{\rho \left( \rho +1\right) }+\sum_{\substack{ \text{Re}%
(\rho )=\frac{1}{2}  \\ Y<\left\vert \gamma \right\vert \leq T}}\frac{\left(
x+h\right) ^{\rho +1}-x^{\rho +1}}{\rho \left( \rho +1\right) }\text{.}
\end{equation*}%
The trivial bound for the first sum on the right hand side is given by%
\begin{equation}
\frac{1}{h}\left\vert \sum_{\substack{ \text{Re}(\rho )=\frac{1}{2}  \\ %
\left\vert \gamma \right\vert \leq Y}}\frac{\left( x+h\right) ^{\rho
+1}-x^{\rho +1}}{\rho \left( \rho +1\right) }\right\vert =O\left( x^{\frac{1%
}{2}}\sum_{\substack{ \text{Re}(\rho )=\frac{1}{2}  \\ \left\vert \gamma
\right\vert \leq Y}}\frac{1}{\left\vert \rho \right\vert }\right) =O\left(
x^{\frac{1}{2}}Y\right) \text{.}  \label{f2}
\end{equation}%
The second sum is split into%
\begin{equation*}
\sum_{\substack{ \text{Re}(\rho )=\frac{1}{2}  \\ Y<\left\vert \gamma
\right\vert \leq T}}\frac{\left( x+h\right) ^{\rho +1}}{\rho \left( \rho
+1\right) }-\sum_{\substack{ \text{Re}(\rho )=\frac{1}{2}  \\ Y<\left\vert
\gamma \right\vert \leq T}}\frac{x^{\rho +1}}{\rho \left( \rho +1\right) }%
\text{.}
\end{equation*}%
Let
\begin{equation*}
D_{Y}^{T}=\left\{ x\in \left[ T,eT\right) :\left\vert \underset{_{\substack{
\text{Re}(\rho )=\frac{1}{2}  \\ Y<\left\vert \gamma \right\vert \leq T}}}{%
\sum }\frac{x^{\rho +1}}{\rho \left( \rho +1\right) }\right\vert >x^{\alpha
}(\log x)^{\beta }(\log \log x)^{\beta }\right\} \text{, }1<\alpha <\frac{3}{%
2}\text{, }\beta >0\text{.}
\end{equation*}%
Then,
\begin{eqnarray*}
\mu ^{\times }D_{Y}^{T} &=&\underset{D_{Y}^{T}}{\int }\frac{dt}{t}=\underset{%
D_{Y}^{T}}{\int }t^{2\alpha }(\log t)^{2\beta }(\log \log t)^{2\beta }\frac{%
dt}{t^{1+2\alpha }(\log t)^{2\beta }(\log \log t)^{2\beta }} \\
&\leq &\frac{1}{(\log T)^{2\beta }(\log \log T)^{2\beta }}\underset{D_{Y}^{T}%
}{\int }\left\vert \sum_{\substack{ \text{Re}(\rho )=\frac{1}{2}  \\ %
Y<\left\vert \gamma \right\vert \leq T}}\frac{t^{\rho +1}}{\rho \left( \rho
+1\right) }\right\vert ^{2}\frac{t^{3-2\alpha }}{t^{4}}dt \\
&=&O\left( \frac{T^{3-2\alpha }}{(\log T)^{2\beta }(\log \log T)^{2\beta }}%
\right) \overset{eT}{\underset{T}{\int }}\left\vert \sum_{\substack{ \text{Re%
}(\rho )=\frac{1}{2}  \\ Y<\left\vert \gamma \right\vert \leq T}}\frac{%
t^{\rho +1}}{\rho \left( \rho +1\right) }\right\vert ^{2}\frac{dt}{t^{4}}%
\text{.}
\end{eqnarray*}%
According to Koyama \cite[p. 79]{K},
\begin{equation*}
\overset{eT}{\underset{T}{\int }}\left\vert \sum_{\substack{ \text{Re}(\rho
)=\frac{1}{2}  \\ Y<\left\vert \gamma \right\vert \leq T}}\frac{t^{\rho +1}}{%
\rho \left( \rho +1\right) }\right\vert ^{2}\frac{dt}{t^{4}}=O\left( \frac{1%
}{Y}\right) \text{.}
\end{equation*}
Taking $Y=T^{3-2\alpha }(\log T)^{1-2\beta }(\log \log T)^{1-2\beta
+\varepsilon }$, we obtain%
\begin{equation*}
\mu ^{\times }D_{Y}^{T}\ll \frac{1}{\log T\left( \log \log T\right)
^{1+\varepsilon }}\text{.}
\end{equation*}
For $n=\left\lfloor \log x\right\rfloor $, $T=e^{n}$, denote $%
E_{n}=D_{Y}^{T} $. Then, $\mu ^{\times }E_{n}\ll \frac{1}{n\left( \log
n\right) ^{1+\varepsilon }}$ and $\mu ^{\times }\cup E_{n}\ll\sum \frac{1}{n\left( \log n\right)
^{1+\varepsilon }}<\infty $.

If $x\in \left[ e^{n},e^{n+1}\right) \setminus E_{n}$, one gets
\begin{equation}
\left\vert \underset{_{\substack{ \text{Re}(\rho )=\frac{1}{2}  \\ %
Y<\left\vert \gamma \right\vert \leq T}}}{\sum }\frac{x^{\rho +1}}{\rho
\left( \rho +1\right) }\right\vert \leq x^{\alpha }(\log x)^{\beta }(\log
\log x)^{\beta }\text{.}  \label{f3}
\end{equation}%
We are interested in achieving $h<x^{\frac{3}{4}}$. If it happens that $x+h\in \left[
e^{n},e^{n+1}\right)\setminus E_{n} $, one shall have
\begin{equation*}
\underset{_{\substack{ \text{Re}(\rho )=\frac{1}{2}  \\ %
Y<\left\vert \gamma \right\vert \leq T}}}{\sum }\frac{(x+h)^{\rho +1}}{\rho
\left( \rho +1\right) }=O\left( x^{\alpha }(\log x)^{\beta
}(\log \log x)^{\beta }\right)
\end{equation*}%
as well, since $x+h<2x$.

The other possibility is that $x+h\in \left[ e^{n+1},e^{n+2}\right) $. In
that case, we proceed as follows%
\begin{equation*}
\underset{_{\substack{ \text{Re}(\rho )=\frac{1}{2}  \\ Y<\left\vert \gamma
\right\vert \leq T}}}{\sum }\frac{(x+h)^{\rho +1}}{\rho \left( \rho
+1\right) }=\underset{_{\substack{ \text{Re}(\rho )=\frac{1}{2}  \\ %
Y<\left\vert \gamma \right\vert \leq eT}}}{\sum }\frac{(x+h)^{\rho +1}}{\rho
\left( \rho +1\right) }-\underset{_{\substack{ \text{Re}(\rho )=\frac{1}{2}
\\ T<\left\vert \gamma \right\vert \leq eT}}}{\sum }\frac{(x+h)^{\rho +1}}{%
\rho \left( \rho +1\right) }\text{.}
\end{equation*}

For $x+h\in \left[ e^{n+1},e^{n+2}\right) \setminus E_{n+1}$, we get%
\begin{equation} \label{f4}
\underset{_{\substack{ \text{Re}(\rho )=\frac{1}{2}  \\ %
Y<\left\vert \gamma \right\vert \leq eT}}}{\sum }\frac{(x+h)^{\rho +1}}{\rho
\left( \rho +1\right) } =O\left( x^{\alpha }(\log x)^{\beta
}(\log \log x)^{\beta }\right) \text{.}
\end{equation}

To estimate $\underset{_{\substack{ \text{Re}(\rho )=\frac{1}{2}  \\ %
T<\left\vert \gamma \right\vert \leq eT}}}{\sum }\frac{(x+h)^{\rho +1}}{\rho
\left( \rho +1\right) }$, we shall consider
\begin{equation*}
D_{T}^{eT}=\left\{ x\in \left[ eT,e^{2}T\right) :\left\vert \underset{
_{\substack{ \text{Re}(\rho )=\frac{1}{2}  \\ T<\left\vert \gamma
\right\vert \leq eT}}}{\sum }\frac{x^{\rho +1}}{\rho \left( \rho +1\right) }%
\right\vert >x^{\alpha }(\log x)^{\beta }(\log \log x)^{\beta }\right\}
\text{.}
\end{equation*}%
By the argumentation above leading to the estimate of $\mu ^{\times
}D_{Y}^{T}$, we get
\begin{equation*}
\mu ^{\times }D_{T}^{eT}\ll \frac{T^{3-2\alpha }}{(\log T)^{2\beta }\left(
\log \log T\right) ^{2\beta }}\cdot \frac{1}{T}\ll \frac{1}{T^{2\alpha -2}}%
\text{.}
\end{equation*}%
Recall that $n=\left\lfloor \log x\right\rfloor $, $T=e^{n}$ and denote $%
F_{n+1}=D_{T}^{eT}$. Notice that $\mu ^{\times }F_{n+1}<\frac{1}{e^{(2\alpha
-2)n}}$ and the series $\sum \frac{1}{e^{(2\alpha -2)n}}$ converges since $%
\alpha >1$. Thus, if we additionally assume $x+h\in \left[
e^{n+1},e^{n+2}\right) \setminus F_{n+1}$, we get%
\begin{equation} \label{f5}
\underset{_{\substack{ \text{Re}(\rho )=\frac{1}{2}  \\ T<\left\vert \gamma
\right\vert \leq eT}}}{\sum }\frac{(x+h)^{\rho +1}}{\rho \left( \rho +1\right) }%
=O\left( x^{\alpha }(\log x)^{\beta }(\log \log x)^{\beta }\right) \text{.}
\end{equation}
Looking back at \eqref{f1}, and taking into account the relations \eqref{f2}, \eqref{f3}, \eqref{f4} and \eqref{f5}, we are left to optimize
\begin{center}
$h$,$\ \frac{x\log x}{h}$, $x^{\frac{1}{2}}Y$ and $\frac{x^{\alpha }(\log
x)^{\beta }(\log \log x)^{\beta }}{h}$,  i.e.,
\end{center}
\begin{center}
$h$, $x^{\frac{1}{2}}\cdot x^{3-2\alpha }(\log x)^{1-2\beta }(\log \log
x)^{1-2\beta }$ and $\frac{x^{\alpha }(\log x)^{\beta }(\log \log x)^{\beta }%
}{h}$
\end{center}
since $\alpha > 1$, $T\approx x$ and $Y=O\left(x^{3-2\alpha }(\log x)^{1-2\beta }(\log \log
x)^{1-2\beta +\varepsilon }\right)$.

Choosing $h\approx x^{\frac{\alpha}{2} }(\log x)^{\frac{\beta}{2} }(\log \log
x)^{\frac{\beta}{2}}$, we get $\frac{1}{2}+3-2\alpha = \frac{\alpha}{2}$ and $1-2\beta=\frac{\beta}{2} $. Hence, $\alpha=\frac{7}{5}$ and $\beta=\frac{2}{5}$. This completes the proof since the sets $E=\cup E_n$ and $F=\cup F_n$ have finite logarithmic measure.

\end{proof}
\end{theorem*}


\begin{thebibliography}{99}
\bibitem{A} M. Avdispahi\'{c}, On Koyama's refinement of the prime geodesic
theorem, arXiv:1701.01642

\bibitem{AG} M. Avdispahi\'{c} and D\v{z}. Gu\v{s}i\'{c}, On the error term
in the prime geodesic theorem, Bull. Korean Math. Soc. \textbf{49} (2012),
no. 2, 367--372.

\bibitem{AS} M. Avdispahi\'{c} and L. Smajlovi\'{c}, An explicit formula and
its application to the Selberg trace formula, Monatsh. Math. \textbf{147}
(2006), no. 3, 183--198.

\bibitem{C} Y. Cai, Prime geodesic theorem, J. Th\'{e}or. Nombres Bordeaux
\textbf{14} (2002), no. 1, 59--72.

\bibitem{D} D.L. DeGeorge, Length spectrum for compact locally symmetric
spaces of strictly negative curvature, Ann. Sci. Ecole Norm. Sup. \textbf{10}
(1977), 133--152.

\bibitem{D1} A. Deitmar, A prime geodesic theorem for higher rank spaces.
Geom Funct Anal. \textbf{14} (2004) 1238--1266

\bibitem{G} P. X. Gallagher, Some consequences of the Riemann hypothesis,
Acta Arith. \textbf{37} (1980), 339--343.

\bibitem{GW} R. Gangolli and G. Warner, Zeta functions of Selberg's type for
some noncompact quotients of symmetric spaces of rank one, Nagoya Math. J.
\textbf{78} (1980), 1--44.

\bibitem{H} D. A. Hejhal, \textit{The Selberg trace formula for $PSL(2,R)$.
Vol I}, Lecture Notes in Mathematics, Vol 548, Springer, Berlin 1976.

\bibitem{H1} H. Huber, Zur analytischen Theorie hyperbolischer Raumformen
und Bewegungsgruppen II, Math. Ann. \textbf{142} (1961), 385--398.

\bibitem{H2} H. Huber, Nachtrag zu \cite{H1}, Math. Ann. \textbf{143}
(1961), 463--464.

\bibitem{I} H. Iwaniec, Prime geodesic theorem, J. Reine Angew. Math. \textbf{349}
(1984), 136--159.

\bibitem{I1} H. Iwaniec, \textit{Spectral Methods of Automorphic Forms (2nd
ed.)}. Amer. Math. Soc., Providence 2002.

\bibitem{K} S. Koyama, Refinement of prime geodesic theorem, Proc. Japan
Acad. Ser A Math. Sci. \textbf{92} (2016), no. 7, 77--81.

\bibitem{LS} W. Luo and P. Sarnak, Quantum ergodicity of eigenfunctions on $%
PSL\_2(Z)\backslash H\symbol{94}2$, Inst. Hautes Etudes Sci. Publ. Math. no.
81 (1995), 207--237.

\bibitem{M} G. A. Margulis, Certain applications of ergodic theory to the
investigation of manifolds of negative curvature, Funkcional. Anal. i
Prilozhen. \textbf{3} (1969), no. 4, 89--90.

\bibitem{P} J. Park, Ruelle zeta function and prime geodesic theorem for
hyperbolic manifolds with cusps, in: G. van Dijk, M. Wakayama (eds.),
\textit{Casimir Force, Casimir Operators and the Riemann Hypothesis}, 9-13
November 2009, Kyushu University, Fukuoka, Japan, Walter de Gruyter 2010.

\bibitem{PP} W. Parry and M. Pollicott, An analogue of the prime number
theorem for closed orbits of Axiom A flows, Ann. of Math. (2) \textbf{118}
(1983), no. 3, 573--591.

\bibitem{PS} M. Pollicott and M. Sharp, Length asymptotics in higher Teichm%
\"{u}ller theory, Proc. Amer. Math. Soc. \textbf{142} (2014), 101--112.

\bibitem{R} B. Randol, On the asymptotic distribution of closed geodesics on
compact Riemann surfaces, Trans. Amer. Math. Soc. \textbf{233} (1977),
241--247.

\bibitem{S} A. Selberg, Harmonic analysis and discontinuous groups in weakly
symmetric Riemannian spaces with applications to Dirichlet series, J. Indian
Math. Soc. B. \textbf{20} (1956), 47--87.

\bibitem{SY} K. Soundararajan and M. P. Young, The prime geodesic theorem,
J. Reine Angew. Math. \textbf{676} (2013), 105--120.
\end{thebibliography}
\end{document}